\documentclass[11pt]{amsart}

\usepackage{amsmath,amssymb,amsthm}
\usepackage{mathrsfs}
\usepackage{graphicx}
\usepackage[margin=1in]{geometry}
\usepackage{macros}
\usepackage[numbers]{natbib}
\usepackage{color}
\usepackage{underscore}
\usepackage{hyperref}
\usepackage{stmaryrd}
\usepackage{tikz}
\usepackage{bussproofs}
\usepackage{bussproofs-extra}
\usepackage{tikz,xcolor,hyperref}
\usepackage{orcidlink}

\newtheorem{theorem}{Theorem}[section]

\newtheorem{lemma}[theorem]{Lemma}
\newtheorem{corollary}[theorem]{Corollary}
\theoremstyle{definition}
\newtheorem{definition}[theorem]{Definition}
\theoremstyle{remark}

\newcommand{\Iota}{\mathrm{I}}

\title{Computational Complexity and Bounded Arithmetic over First-Order Structures}

\author{Shogo Saitou \orcidlink{0009-0004-4157-8076}}
\address{Tohoku University \\ Mathematical Institute}
\email{palalansouki@gmail.com}

\date{\today}
 
\begin{document}

\begin{abstract}
  We study computation over a first-order structure $\mathbb{K}$ using $\mathbb{K}$-alphabetic Turing machines.
  These generalize ordinary Turing machines by using tuples of elements of $\mathbb{K}$ as tape symbols and
  allowing parameter-free quantifier-free definable unit-time operations.
  We prove an analogue of the Cook--Levin theorem,
  identifying the existential theory of $\mathbb{K}$ with parameters as a complete problem for $\NP[\mathbb{K}]$.
  For $\mathbb{K}=\mathbb{R}$, the Boolean restriction of this theorem identifies
  the existential theory of the reals as a complete problem for $\NP[\Real]\boolRestr$.
  
  For a theory $T$ with $\mathbb{K}\models T$, we then introduce a two-sorted theory $T + \NIA{\Gamma}$.
  Under some assumptions on $T$, we prove that the $\nSigma{i}$-definable functions in $T + \NIA{\nSigma{i}}$ are exactly
  the functions computable in polynomial time with an oracle in the $(i-1)$st level of the $\mathbb{K}$-alphabetic polynomial-time hierarchy.
\end{abstract}

\maketitle

\section{Introduction}

Classical complexity theory studies computations on finite strings over $\Bool$,
classified by resource bounds such as time and space.
The important complexity classes are $\ClassP$ and $\NP$,
which are the classes of problems solvable in polynomial time by deterministic and nondeterministic Turing machines.
The satisfiability problem $\SAT$ is intensively studied as a complete problem for $\NP$.

This can be viewed from a more algebraic perspective;
if $\Bool$ is regarded as the two-element Boolean structure, then $\SAT$ is the existential theory of $\Bool$.
This can be extended to an arbitrary first-order structure, say, an $\mathcal{L}$-structure $\mathbb{K}$:
one may consider the set of quantifier-free $\mathcal{L}$-formulas which are satisfiable in $\mathbb{K}$.
For $\mathcal{L} = \mathcal{L}_\RCOF := \{0,1,+,\cdot,- ^{-1},<,=\}$ and $\mathbb{K} = \Real$ this gives the \emph{existential theory of the reals} $\ETh{\Real}$.
By Tarski's theorem the full first-order theory of real closed fields is decidable~\cite{Tarski1931},
while the exact complexity of its existential fragment is still not known.
Many natural geometric and algebraic problems are complete for this existential theory~\cite{MeerWurm2025}.

This analogy between $\SAT$ and $\ETh{\Real}$ motivates the study of algebraic computation models.
Blum, Shub and Smale~\cite{BlumShubSmale1990} introduced a model of computation over the real numbers,
now known as the BSS machine, which serves as a foundation for subsequent studies of algebraic computation models.
They defined the corresponding complexity classes $\ClassP_{\Real}$ and $\NP_{\Real}$,
and proved an analogue of the Cook--Levin theorem for $\NP_{\Real}$,
identifying the existential theory of the reals as a complete problem.
Goode~\cite{Goode1994} generalized the BSS machine to first-order structures,
and defined the corresponding generalized complexity classes, e.g., $\ClassP_{\mathbb{K}}$ and $\NP_{\mathbb{K}}$.
Later, Bournez, Cucker, Jacobé de Naurois, and Marion~\cite{BournezCuckerJacobedeNauroisMarion2005}
gave an implicit characterization of the polytime computable functions over a first-order structure by generalizing the notion of a safe function.

These approaches define computation over a first-order structure $\mathbb{K}$
in terms of a prescribed collection of primitive operations over finite tuples of elements of $\mathbb{K}$.
The machine configuration at each step is represented by such a tuple.
For example, addition, multiplication, and tests for positivity are taken to have unit cost.
This raises a question of canonicity: which operations should count as elementary computations over the structure?
Our approach is just to allow all parameter-free quantifier-free definable operations as unit-time operations.
This reduces the problem of canonicity to a (more manageable) problem of choosing an appropriate language.
Furthermore, this is naturally justified from a proof-theoretic perspective.

We introduce the \emph{$\mathbb{K}$-alphabetic Turing machine},
which allows any transition whose graph is definable by a parameter-free quantifier-free $\mathcal L$-formula
as a unit-time operation.
The definition requires only a slight modification of that of an ordinary Turing machine.
For this model, we define corresponding complexity classes $\ClassP[\mathbb{K}]$, $\NP[\mathbb{K}]$, etc.,
and prove an analogue of the Cook--Levin theorem:
the existential theory of $\mathbb{K}$ with parameters is complete for the nondeterministic polynomial-time class $\NP[\mathbb{K}]$.
When $\mathcal{L} = \mathcal{L}_\Bool := \{0,1,\land,\lor,\lnot\}$ and $\mathbb{K} = \Bool$,
the model coincides with the ordinary Turing-machine model, and the corresponding complexity classes coincide as well.

This model does not allow parameters but allows more unit-time operations than the BSS model.
However, in the parameter-free case, it is easy to see that $\NP_\Real^{\text{pf}} = \NP[\Real]$,
because the enriched set of elementary operations can be simulated using nondeterministic choices.

The definition of the machine is motivated by bounded arithmetic.
Buss~\cite{Buss1986} characterized polynomial-time function classes by the provably total functions of weak fragments of arithmetic.
We introduce a two-sorted language with an alphabet sort, interpreted by $\mathbb{K}$,
and a string sort for finite strings over that alphabet.
The resulting theory $T + \NIA{\Gamma}$ consists of an alphabet theory $T$
and induction on notation for formulas in $\Gamma$.
Under suitable hypotheses on $T$, we show that, for every $i\geq 1$,
the provably total functions of $T+\NIA{\nSigma{i}}$ are exactly
the functions in $\Square{i}[\mathbb{K}]=\FP[\mathbb{K}]^{\ClassSigma{i-1}[\mathbb{K}]}$.
In particular, this gives a proof-theoretic characterization of these function classes over the real numbers by taking $T=\RCOF$.

\section{Computation over First-Order Structures}

Let $\mathcal L$ be a first-order finite language containing constants $0$ and $1$,
and let $\mathbb{K}$ be a $\mathcal L$-structure in which $0^{\mathbb{K}} \neq 1^{\mathbb{K}}$.
The two basic examples used throughout the paper are the two-element Boolean algebra and the field of real numbers.
$\Bool = \{0, 1\}$ is the structure of two elements in the language of Boolean algebras
$\mathcal L_{\BOOL} = \{0,1,-,+,\mathord{\cdot},=\}$.
$\Real$ is the structure of real numbers in the language of ordered fields $\mathcal L_{\RCOF} = \{0, 1, \cdot, +, -, ^{-1}, <, =\}$.

\begin{definition}[Existential Theory of a Structure]
  The \emph{theory} and \emph{existential theory} of $\mathbb{K}$ are, respectively, the sets
  \begin{align*}
    \Th{\mathbb{K}} &= \{A : \text{$\mathcal L$-sentences} \mid \mathbb{K} \models A\},\\
    \ETh{\mathbb{K}} &= \{A[\vec{x}] : \text{quantifier-free $\mathcal L$-formula} \mid \mathbb{K} \models \exists\vec{x} A[\vec{x}]\}.
  \end{align*}
\end{definition}

Thus the satisfiability problem $\SAT$ can be identified with $\ETh{\Bool}$.
Its real counterpart $\ETh{\Real}$ is the \emph{existential theory of the reals}.

Tarski proved that $\Th{\Real}$ is decidable~\cite{Tarski1931}.
Nevertheless, the exact complexity of the existential theory of the reals $\ETh{\Real}$ remains unknown;
it is $\NP$-hard and belongs to $\PSPACE$, as shown by Canny~\cite{Canny1988}.
Many geometric and algebraic problems complete for $\ETh{\Real}$ are known. See~\cite{MeerWurm2025} for a survey.
This analogy between $\SAT$ and $\ETh{\Real}$ motivates a uniform model of computation over first-order structures,
which we introduce in the next section.

\subsection{$\mathbb{K}$-alphabetic Turing machines}

For a $\mathcal L$-structure $\mathbb{K}$, a $\mathbb{K}$-alphabetic Turing machine is defined in the same way as an ordinary Turing machine,
except that its tape symbols and states are now tuples of $\mathbb{K}$ and its transition function is quantifier-free $\mathcal L$-definable.
For the technical details of ordinary Turing machines, see~\cite{AroraBarak2009}.

\begin{definition}[$\mathbb{K}$-alphabetic Turing Machine]
  Fix $k \ge 1$.  A \emph{$\mathbb{K}$-alphabetic $k$-tape Turing machine},
  or simply \emph{$\mathbb{K}$-alphabetic TM}, is a triple $M = (\Gamma,Q,\delta)$ such that
  \begin{enumerate}
    \item $\Gamma = \mathbb{K}^g$ is the tape alphabet for some $g \ge 1$;
    \item $Q = \mathbb{K}^q$ is the set of states for some $q \ge 1$, with distinguished
      and distinct states $\mathtt{start}$ and $\mathtt{halt}$ in $Q$, both coded by tuples built from $0^\mathbb{K}$ and $1^\mathbb{K}$;
    \item The transition function
      \[
        \delta: \Gamma^k \times Q \to
        \Gamma^{k-1} \times Q \times \{\mathtt{left},\mathtt{stay},\mathtt{right}\}^k
      \]
      has a graph defined by a parameter-free quantifier-free $\mathcal L$-formula.
      The three head directions are coded by distinct tuples built from $0^\mathbb{K}$ and $1^\mathbb{K}$.
  \end{enumerate}
\end{definition}

The first tape is read-only, and the remaining tapes are working tapes. All heads initially point to the leftmost cell,
and the working tapes are blank.
At each computation step, the machine reads the symbols under the $k$ heads and the current state, updates the symbols under the $k-1$ working heads, moves all $k$ heads, and updates the state according to the transition function $\delta$.
The computation halts upon entering state $\mathtt{halt}$; its output is the string on the first working tape.
Computation time is the number of steps before halting.
If $\mathbb{K} = \Bool$, then a $\mathbb{K}$-alphabetic TM is just an ordinary Turing machine.

Oracle computation is defined as usual: an additional oracle tape and a query state allow the machine state to branch
according to whether the string on the oracle tape belongs to a given oracle language.

A function $f\colon \mathbb{K}^{<\omega} \times \dots \times \mathbb{K}^{<\omega} \to \mathbb{K}^{<\omega}$ is \emph{polynomial-time computable}
if a $\mathbb{K}$-alphabetic TM $M$ computes it in time polynomial in the total lengths of its inputs.
We denote the resulting class by $\FP[\mathbb{K}]$, and its relativization to an oracle $\mathcal O$ by $\FP[\mathbb{K}]^{\mathcal O}$.

The truth evaluation of a quantifier-free $\mathcal L$-formula with parameters from $\mathbb{K}$ can be computed in linear time by a $\mathbb{K}$-alphabetic TM;
the algorithm is a straightforward evaluation of the Reverse Polish notation of formulas.
This fact will be used later.

\subsection{Complexity classes}
A $\mathbb{K}$-alphabetic complexity class is a class of subsets of $\mathbb{K}^{<\omega}$.
For an ordinary complexity class $\mathbf{C}$, we denote its $\mathbb{K}$-alphabetic version by $\mathbf{C}[\mathbb{K}]$.
It is easy to check that $\mathbf{C}[\Bool] = \mathbf{C}$ for any standard complexity class $\mathbf{C}$.

We define the class $\ClassP[\mathbb{K}] = \ClassSigma{0}[\mathbb{K}] = \ClassPi{0}[\mathbb{K}]$ as the class of sets
decided in polynomial time by a $\mathbb{K}$-alphabetic TM.
For $i \ge 1$, a set $A$ belongs to $\ClassSigma{i}[\mathbb{K}]$ if there are a polynomial $p$ and a polynomial-time function $f$ such that

\begin{align*}
  x \in A \iff \exists y_1 \in \mathbb{K}^{p(|x|)} \forall y_2 \in \mathbb{K}^{p(|x|)} \dots Q_i y_i \in \mathbb{K}^{p(|x|)} \, [f(x,y_1,\dots,y_i) = 1]
\end{align*}

$Q_i = \exists$ if $i$ is odd, and $Q_i = \forall$ if $i$ is even. The class $\ClassPi{i}[\mathbb{K}]$ is defined similarly, with the first quantifier being universal.
We put

\begin{align*}
  \NP[\mathbb{K}] := \ClassSigma{1}[\mathbb{K}] \quad \quad \quad \quad
  \Square{i}[\mathbb{K}] := \FP[\mathbb{K}]^{\ClassSigma{i-1}[\mathbb{K}]} \quad (i \ge 1)
\end{align*}

If $\mathbf C$ is a complexity class over $\mathbb{K}$, its restriction to Boolean inputs is

\begin{align*}
  \mathbf C\boolRestr := \{A \cap \Bool^{<\omega} : A \in \mathbf C\}
\end{align*}

Here $\Bool$ is identified with $\{0^{\mathbb{K}}, 1^{\mathbb{K}}\}$.
This restriction lets us compare a problem with a Boolean encoding, such as an existential theory, to $\mathbb{K}$-alphabetic classes.

\subsection{Cook--Levin theorem}

To state the analogue of the Cook--Levin theorem, we distinguish existential theories with parameters and without parameters.

\begin{align*}
  \EXSp{\mathbb{K}} &:= \{F[\vec x] : \text{quantifier-free $\mathcal L$-formula with parameters from $\mathbb{K}$} \mid \mathbb{K} \models \exists \vec{x} \, F[\vec x]\} \\
  \EXS{\mathbb{K}} &:= \EXSp{\mathbb{K}} \cap \Bool^{<\omega}
\end{align*}

The linear-time evaluation described above shows that $\EXSp{\mathbb{K}} \in \NP[\mathbb{K}]$
and that $\EXS{\mathbb{K}} \in \NP[\mathbb{K}] \boolRestr$.
The next theorem is the analogue of the Cook--Levin theorem for our computation model.
It is stated in both parameterized and unparameterized versions.

\begin{theorem}[Cook--Levin theorem over $\mathbb{K}$]
  \leavevmode
  \begin{enumerate}
    \item $\EXSp{\mathbb{K}}$ is complete for $\NP[\mathbb{K}]$ under polynomial-time computable functions over $\mathbb{K}$.
    \item $\EXS{\mathbb{K}}$ is complete for $\NP[\mathbb{K}] \boolRestr$ under polynomial-time computable functions over $\mathbb{K}$.
  \end{enumerate}
\end{theorem}
\begin{proof}[Proof sketch]
  We will show (1). (2) is proved similarly.
  The proof is almost identical to that of the ordinary Cook--Levin theorem; see~\cite{AroraBarak2009} for details.
  
  We prove hardness by the usual Cook--Levin construction, with the transition predicates interpreted in $\mathbb{K}$.
  Let $A \in \NP[\mathbb{K}]$.
  Choose a polynomial-time $\mathbb{K}$-alphabetic TM $M$ and a polynomial $p$ such that
  \[
    x \in A
    \quad\Longleftrightarrow\quad
    \exists y \in \mathbb{K}^{p(\size{x})}\,[M(x,y)=1].
  \]
  $M$ is a $3$-tape machine with two input tapes and one working tape.
  Moreover, we assume that $M$ is oblivious, that is,
  the head positions at time $i$ are computable from only the length of inputs and $i$ independently of the contents of the tapes~\cite{AroraBarak2009}.
  Let $T(n)$ be a polynomial-time bound for $M$, with $n$ being the length of the input $x$.

  We define an open $\mathcal L$-formula which describes snapshots of the computation of $M$.
  A \emph{snapshot} is a pair consisting of the current state and the symbol under the head of the working tape.
  We use the following open $\mathcal L$-formulas:
  \[
  \begin{tabular}{@{}l p{0.66\linewidth}@{}}
    $\Tr[z_1,z_0,z_{-1},\gamma_0,\gamma_1]$
      & states that $z_1$ is the next snapshot of the computation from the current snapshot $z_0$,
        where $\gamma_0$ and $\gamma_1$ are the symbols under the two input-tape heads,
        and $z_{-1}$ is the snapshot at the last visit to the present work-tape cell; \\[1mm]
    $\mathrm{Start}[z]$
      & states that $z$ is the initial snapshot; \\
    $\mathrm{Halt}[z]$
      & states that $z$ is the halting snapshot.
  \end{tabular}
  \]
  In the definition of $\Tr$, if the work-tape head has not visited the present cell before, we put $z_{-1}=z_0$.
  These formulas are open because the transition function of $M$ is defined by an open $\mathcal L$-formula.

  For each $x \in \mathbb{K}^{<\omega}$, define the open $\mathcal L$-formula
  \[
  \begin{aligned}
    \Phi_x[y,z_0,\ldots,z_{T(\size{x})}]
      :=&\quad
      \mathrm{Start}[z_0] \land \mathrm{Halt}[z_{T(\size{x})}] \\
      &{}\land
      \bigwedge_{i<T(\size{x})}
        \Tr[
          z_{i+1},
          z_i,
          z_{r(\size{x},i)},
          (x)_{h_0(\size{x},i)},
          (y)_{h_1(\size{x},i)}
        ].
  \end{aligned}
  \]
  Here $h_k(n,i)$ is the position of the head of tape $k \in \{0,1,2\}$
  (tape $0$ for input $x$; tape $1$ for input $y$; tape $2$ for the work tape) at step $i$ on inputs of length $n$,
  and $r(n,i)$ is the last step before $i$ at which the head of the working tape visited the position $h_2(n,i)$.
  Since $M$ is oblivious and $T(n)$ is polynomially bounded,
  the map $x \mapsto \Phi_x$ is computable in polynomial time by a $\mathbb{K}$-alphabetic TM.
  Moreover,
  \[
  \begin{aligned}
    x \in A
    &\iff
      \exists y \in \mathbb{K}^{p(\size{x})}\,[M(x,y)=1] \\
    &\iff
      \mathbb{K} \models
      \exists y\,\exists z_0\cdots\exists z_{T(\size{x})}\,
      \Phi_x[y,z_0,\ldots,z_{T(\size{x})}] \\
    &\iff
      \Phi_x[y,z_0,\ldots,z_{T(\size{x})}] \in \EXSp{\mathbb{K}}.
  \end{aligned}
  \]
  Thus $A$ is reducible to $\EXSp{\mathbb{K}}$ by a polynomial-time computable
  $\mathbb{K}$-function.  If $x \in \Bool^{<\omega}$, then the same construction
  has no parameters other than the constants $0^{\mathbb{K}}$ and
  $1^{\mathbb{K}}$, and hence gives the reduction to $\EXS{\mathbb{K}}$.
\end{proof}

\section{Bounded Arithmetic over a First-Order Structure}

Bounded arithmetic is often presented as a collection of fragments of first-order arithmetic. 
Its two-sorted presentation instead regards numbers and their finite sets as separate objects; see~\cite{Zambella1996} and~\cite{CookNguyen2010}.
For the present purpose even the number sort is unnecessary.
Following the \emph{induction-on-notation} approach of Ferreira~\cite{Ferreira1994},
we use one sort (\emph{alphabet}) for elements of an arbitrary first-order structure
and another sort (\emph{string}) for finite strings of such elements.
When the alphabet is $\Bool = \{0,1\}$, the construction specializes to the familiar theory of finite binary strings.

Fix a one-sorted language $\mathcal L$ with a distinguished constant $c$,
and let $\mathcal N(\mathcal L)$ be the two-sorted language whose sorts are called \emph{alphabet} and \emph{string}.
Lower-case variables $x,y,z,...$ range over the alphabet sort and upper-case variables $X,Y,Z,...$ over the string sort.
The alphabet-sort $\mathcal L$-symbols are precisely those of $\mathcal L$.

On strings we have the empty string $\emptyseq$, the singleton operation $\brak{x}$, concatenation $X^\frown Y$, repetition $X \times Y$,
and the prefix relation $X \subseteq Y$.
Here $X \times Y$ denotes $\size{Y}$ copies of $X$.
We usually omit singleton brackets and concatenation symbols.

\subsection{Theory}

The basic theory of strings is universal and consists of the following axioms

\begin{gather*}
  X \emptyseq = X \quad \quad
  X (Y z) = (X Y) z \quad \quad
  X \times \emptyseq = \emptyseq \quad \quad
  X \times (Y z) = (X \times Y) X \\
  X \subseteq \emptyseq \liff X = \emptyseq \quad \quad
  X \subseteq Y z \liff X \subseteq Y \lor X = Y z \\
  X z \ne \emptyseq \quad \quad
  X x = Y y \liff X = Y \land x = y
\end{gather*}

These symbols and axioms allow the usual elementary operations on strings to be defined without a length (\emph{number}) sort.
In particular, the relation defined as

\begin{align*}
  X \le Y := \brak{c} \times X \subseteq \brak{c} \times Y
\end{align*}

represents the usual length comparison, and we write $\exists X \le T\, A$ and $\forall X \le T\, A$
for the \emph{bounded quantification} $\exists X (X \le T \land A)$ and $\forall X (X \le T \to A)$.
We also define the subword relation by

\begin{align*}
  X \lsubseteq Y := \exists Z\, (Z \subseteq Y \land ZX \subseteq Y)
\end{align*}

and we write $\exists X \lsubseteq T\, A$ and $\forall X \lsubseteq T\, A$
for the \emph{sharply bounded quantification} $\exists X (X \lsubseteq T \land A)$ and $\forall X (X \lsubseteq T \to A)$.

We work with formulas in negation-normal form.
The class $\nSigma{0} = \nPi{0}$ is defined as the smallest set of formulas containing atomic formulas
and closed under connectives and sharply bounded quantification.
Inductively, $\nSigma{i+1}$ (resp. $\nPi{i+1}$) is the smallest set of formulas containing $\nSigma{i} \cup \nPi{i}$
and closed under connectives, sharply bounded quantification, and existential (resp. universal) bounded quantification.
A formula is $\nDelta{i}$ over a theory $T$ if it is provably equivalent in $T$ to both a $\nSigma{i}$-formula and a $\nPi{i}$-formula.

For a class $\Gamma$ of two-sorted formulas, \emph{induction on notation for $\Gamma$} is the scheme

\begin{align*}
  A[\emptyseq] \land \forall X \forall x\, (A[X] \to A[X x]) \to \forall X\, A[X]
\end{align*}

for all $A[X] \in \Gamma$ ($A$ may contain parameters).
We write $\NIA{\Gamma}$ for the basic theory of strings which is stated above together with this schema.
For a $\mathcal L$-theory $T$, we denote by $T + \NIA{\Gamma}$ the two-sorted theory obtained by adding $\NIA{\Gamma}$ to $T$,
assuming that $T$ is formulated in the alphabet sort.

\subsection{Two Basic Alphabet Theories}

We now specify the two important examples.
They recover ordinary binary computation and computation over the real numbers, respectively.

\begin{definition}[The theory $\BOOL$]
  In the language
  \[
    \mathcal L_{\BOOL}=\{0,1,-,+,\mathord{\cdot},=\},
  \]
  the theory $\BOOL$ consists of the equality axioms and the universal
  closures of
  \begin{gather*}
    x=0\lor x=1,\qquad 0\ne1,\qquad
    -0=1,\qquad -1=0,\\
    x\cdot0=0,\qquad x\cdot1=x,\qquad
    x+0=x,\qquad x+1=1.
  \end{gather*}
\end{definition}

The theory $\BOOL$ has, up to isomorphism, the two-element Boolean algebra $\Bool = \{0, 1\}$ as its unique model.
The theory $\BOOL+\NIA{\Gamma}$ is equivalent to the binary string theory of Ferreira~\cite{Ferreira1994}.

\begin{definition}[The theory $\RCOF$]
  Let
  \[
    \mathcal L_{\RCOF}
      =\{0,1,+,\mathord{\cdot},-,{}^{-1},<,=\}.
  \]
  The theory of \emph{real closed ordered fields} $\RCOF$ consists of the universal closures
  of
  \begin{gather*}
    (x+y)+z=x+(y+z),\qquad x+y=y+x,\qquad
    x+0=x,\qquad x+(-x)=0,\\
    (x\cdot y)\cdot z=x\cdot(y\cdot z),\qquad
      x\cdot y=y\cdot x,\qquad
    x\cdot1=x,\qquad x\ne0\to x\cdot x^{-1}=1,\\
    x\cdot(y+z)=x\cdot y+x\cdot z,\qquad 0^{-1}=0,\qquad
    0\ne1,\qquad \neg(x<x),\\
    (x<y\land y<z)\to x<z,\qquad
      x<y\lor x=y\lor y<x,\\
    x<y\to x+z<y+z,\qquad
      0<x\land0<y\to0<x\cdot y, \\
    \exists y\,(x=y\cdot y\lor -x=y\cdot y),\\
    \exists y\,
      \left(y^{2n+1}+x_{2n}y^{2n}+\cdots+x_1y+x_0=0\right)
      \quad(\text{for all $n \in \mathbb N$}).
  \end{gather*}
\end{definition}

In particular, $\Real\models\RCOF$, and $\RCOF$ is complete and admits quantifier elimination~\cite{Tarski1931}.
$\BOOL$ and $\RCOF$ provide, respectively, the binary and real-closed-field alphabets used below.

For now, we write $\Alp$ for the set of underlying elements described by $T$.
Thus $\Square{i}[\Alp]$ is the internalized version of $\Square{i}[\mathbb{K}]$, with $\mathbb{K}$ being the alphabet $\Alp$,
which will be discussed in next sections.

\subsection{Internalizing Alphabetic Turing Machines}

In this subsection we internalize $\Alp$-alphabetic Turing-machine computations in the theory $T + \NIA{\nSigma{1}}$.
Throughout this subsection we assume that $\mathcal L$ contains constants $0$ and $1$, and that $T \vdash 0 \neq 1$ (and we let $c := 1$).

We write $X<Y$ for $X\le Y \land \neg(Y\le X)$ and $X \equiv Y$ for $X \le Y \land Y \le X$.
We use block coding for sequences and the associated projection function.
The notation $(X)^L_I=Y$ means that $Y$ is the $\size{I}$-th block of $X$ when $X$ is divided into blocks of length $\size{L}$, and abbreviates the following formula

\begin{align*}
  (I1 \times L \le X \land\ &\ 
    \exists X' \subseteq X\, (X' \equiv I \times L \land Y \equiv L \land X'Y \subseteq X)) \\
  {}\lor\ (X < I1 \times L \ \land\ & Y=\emptyseq)
\end{align*}

\begin{lemma}
  The function $(X,L,I) \mapsto (X)^L_I$ is provably total in
  $\NIA{\nSigma{1}}$.
\end{lemma}

For simplicity of notation, we state the construction for one-tape machines, with tape alphabet $\Gamma = \Alp$.
Extending the construction to general $k$-tape machines is straightforward and omitted.
Let $M$ be an $\Alp$-alphabetic TM,
whose transition graph is defined by a parameter-free quantifier-free $\mathcal L$-formula, and suppose that the totality and
uniqueness of values of this transition relation are provable in $T$.
Let $Q$ be a fixed unary string coding the block length used for states.

We define a $\nSigma{0}$-formula

\begin{align*}
  \Comp_M(S,T,X,W,H) :=\ &
    \InitTape(X,(W)^S_{\emptyseq})
    \land \InitState((H)^Q_{\emptyseq}) \\
  {} \land\,&
    \forall I \subsetneq T\,
      \Tr \left((W)^S_I,(H)^Q_I,(W)^S_{I1},(H)^Q_{I1}\right)
\end{align*}

expressing that $W$ and $H$ code the tape configurations and states of the run of $M$ on input $X$ through time $T$,  restricted to space $S$
(for time $T$ and space $S$, we are only interested in their lengths).
The intended reading is that $(W)^S_I$ is the tape configuration at time $\size{I}$, and $(H)^Q_I$ is the state at time $\size{I}$.

Here $\InitTape(X,V)$ says that $V$ is the initial tape configuration with input $X$,
$\InitState(U)$ says that $U$ is the initial state,
and $\Tr(V_0,U_0,V_1,U_1)$ says that $(V_1,U_1)$ is obtained from $(V_0,U_0)$ by one transition of $M$ if the state is not $\mathtt{halt}$,
and that $(V_1,U_1)=(V_0,U_0)$ if the state is $\mathtt{halt}$.
Since the transition function of $M$ is given by a quantifier-free $\mathcal L$-formula, these auxiliary formulas may be chosen in $\nSigma{0}$.

Let $\OutTape(Y,V)$ be a $\nSigma{0}$-formula expressing that $Y$ is the output tape of a configuration $V$.
We define the result of computation $Y$ with input $X$ (in space $S$, at time $T$) by

\begin{align*}
  \CompRes_M(S,T,X,Y) :=\ &
    \exists W \le S \times T1\, \exists H \le Q \times T1\, \\
      &\left(\Comp_M(S,T,X,W,H) \land \OutTape(Y, (W)^S_T) \right)
\end{align*}

Thus $\CompRes_M$ is a $\nSigma{1}$-formula.

The oracle computation is internalized in the same way.
If $\Phi[Z]$ is a two-sorted formula with one free string variable, then $\Comp^{\Phi}_M(S,T,X,W,H)$ is
obtained from $\Comp_M$ by replacing the transition clause with following.

\begin{align*}
  \forall I \subsetneq T\,\bigl[
    (\Phi[(W)^S_I] \to\ &
      \Tr_1((W)^S_I,(H)^Q_I,(W)^S_{I1},(H)^Q_{I1})) \\
  {}\land
    (\neg\Phi[(W)^S_I] \to\ &
      \Tr_0((W)^S_I,(H)^Q_I,(W)^S_{I1},(H)^Q_{I1}))\bigr]
\end{align*}

The corresponding result predicate is similarly defined as

\begin{align*}
  \CompRes^{\Phi}_M(S,T,X,Y) :=\ &
    \exists W \le S \times T1\, \exists H \le Q \times T1\, \\
      &\left(\Comp^{\Phi}_M(S,T,X,W,H) \land \OutTape(Y, (W)^S_T) \right)
\end{align*}

If $\Phi$ is a $\nSigma{i}$-formula, then $\CompRes^\Phi_M$ is a $\nSigma{i+1}$-formula.

In the next theorem, the notation ${X}^{\times n}$ denotes the $n$-fold
iteration of the $\times$-operation: $X^{\times n} := \underbrace{X \times \cdots \times X}_n$.

\begin{theorem}[Internalization of $\mathbb{K}$-alphabetic computation]
  Assume that $\mathbb{K}$ is a model of $T$ and that $T$ is $\forall\exists$-complete for $\mathbb{K}$, that is,
  for every open formula $A[\vec x, \vec y]$, $\mathbb{K} \models \forall \vec x\,\exists \vec y\,A[\vec x, \vec y]$ implies $T \vdash \forall \vec x\,\exists \vec y\,A[\vec x, \vec y]$.
  Then
  \begin{enumerate}
    \item For every $\mathbb{K}$-alphabetic TM $M$ run in time $n^k + c$,
      \[
        T+\NIA{\nSigma{1}}
        \vdash \forall X\,\exists Y\,
          \CompRes_M({X}^{\times k}1^c,{X}^{\times k}1^c,X,Y)
      \]
    \item Let $\Phi[X]$ be a $\nSigma{i}$-formula, and let $M$ be a $\mathbb{K}$-alphabetic TM with oracle $\Phi$ run in time $n^k + c$.  Then
      \[
        T+\NIA{\nSigma{i+1}}
        \vdash \forall X\,\exists Y\,
          \CompRes^\Phi_M({X}^{\times k}1^c,{X}^{\times k}1^c,X,Y)
      \]
  \end{enumerate}
\end{theorem}

\begin{proof}
  We prove (1).  As in the definition of $\CompRes_M$, it suffices to prove
  \[
    T+\NIA{\nSigma{1}}
    \vdash
    \forall X\,\exists W\,\exists H\,
      \Comp_M(X^{\times k}1^c,{X}^{\times k}1^c,X,W,H),
  \]
  Let us denote $S, T := X^{\times k}1^c$.
  Fix $X$.  We show this by induction on notation for the
  $\nSigma{1}$-formula
  \[
    B[I] :=
    I\le T \to \exists W \le S\times I1\,\exists H \le Q\times I1\, \Comp_M(S,I,X,W,H).
  \]
  The base case $B[\emptyseq]$ follows from the definitions of $\InitTape$ and $\InitState$.
  For the induction step, assume $B[I]$, and by the induction hypothesis,
  we have $W_0$ and $H_0$ such that $\Comp_M(S,I,X,W_0,H_0)$ holds.
  Since the totality of the transition function is provable in $T$, there exist blocks $V$ and $U$
  such that $\Tr\left((W_0)^S_I,(H_0)^Q_I,V,U\right)$ holds.
  Define $W_1:=W_0V$ and $H_1:=H_0U$.  Then $\Comp_M(S,Iz,X,W_1,H_1)$ holds by definition.

  The proof of (2) is similar.  We use induction on notation for
  \[
    B[I] :=
    I\le T \to \exists W \le S\times I1\,\exists H \le Q\times I1\,\Comp^\Phi_M(S,I,X,W,H).
  \]
  This formula is in $\nSigma{i+1}$.  The induction step is the same as above,
  except that the $V$ and $U$ are chosen according to $\Phi[(W_0)^S_I]$.
\end{proof}

Consequently, every alphabetic polynomial-time computation (with oracle) can be internalized in the corresponding theory. 
In the presence of an oracle at level $i-1$, this yields provable totality in the $i$-th fragment.

\begin{corollary}\label{col:internalization}
  Assume that $T$ is $\forall\exists$-complete for $\mathbb{K}$.
  For every $i \ge 1$, every function in $\Square{i}[\mathbb{K}]$ is provably total in $T+\NIA{\nSigma{i}}$.
\end{corollary}

Moreover, we allow an $\Alp$-alphabetic TM to take $\Alp$-elements as inputs and outputs
by representing each element on a separate tape of length $1$.
This gives a $\mathbb{K}$-alphabetic TM computable function in
$\mathbb{K} \times \cdots \times \mathbb{K} \to \mathbb{K} \times \cdots \times \mathbb{K}$.

\section{Provably Total Functions of \texorpdfstring{$T+\NIA{\nSigma{i}}$}{T+NIA(Sigma i)}}

We will show the converse of Corollary~\ref{col:internalization},
that is, every function provably total in $T+\NIA{\nSigma{i}}$ belongs to $\Square{i}[\Alp]$.
First we introduce the notion of a $\forall\exists!$-axiomatizable theory, which is needed for the witnessing theorem.

\begin{definition}[$\forall\exists!$-axiomatizable theory]
  A theory $T$ is \emph{$\forall\exists!$-axiomatizable}
  if $T$ can be axiomatized by a set of universal sentences and sentences of the form
  $\forall \vec x\,\exists \vec y\,A[\vec x,\vec y]$, where $A$ is
  quantifier-free, and for every such sentence
  $T$ proves the uniqueness of $\vec y$ for each $\vec x$, that is,
  \[
    T \vdash \forall \vec x\,\forall \vec y_1\,\forall \vec y_2\, (A[\vec x,\vec y_1] \land A[\vec x,\vec y_2] \to \vec y_1 = \vec y_2).
  \] 
\end{definition}

The theories $\BOOL$ and $\RCOF$ are both
$\forall\exists!$-axiomatizable.
For $\RCOF$, uniqueness is obtained by choosing the least root in the field order
whenever an axiom asserts the existence of a root.

\subsection{Witnessing theorem}

We use the one-sided sequent calculus for formulas in negation-normal form.
Let $\operatorname{pair}(X,Y)$ be a polynomial-time pairing function for strings,
with projections $\pi_0$ and $\pi_1$.
$\operatorname{nth}(I,Y)$ denotes the $I$-th member of the coded string-indexed sequence of strings $Y$.
Finally, $\operatorname{head}(X)$ denotes the first alphabet element occurring in $X$.
We let all of these functions return predetermined default values (e.g. $0$ or $\emptyseq$) when the input is invalid.
By the previous section, these polynomial-time functions are provably total in $T+\NIA{\nSigma{1}}$.
Furthermore, we assume that these functions use sufficiently efficient encodings;
the size of pairs and coded string-indexed sequences is linear in the sum of the sizes of their components.

We also use notation $\size{X}$ for the length of a string $X$ and (in)equalities between those lengths.
Since our theory does not include number-sort, this does not make sense strictly,
but its meaning is clear; replacing the inequality symbol with $\le$, the length addition with $^\frown$, and the length multiplication with $\times$, etc., 
yields a formal statement in our theory.
For a tuple of strings $\vec{X} = X_1, \ldots, X_n$, write $\size{\vec X}:=\sum_j\size{X_j}$.

\begin{definition}[The classes $\CE{i}$ and $\CA{i}$]
  Let $i\ge 1$.  The positive class $\CE{i}$ (resp. negative class $\CA{i}$) is the smallest class of formulas which
  contains $\nSigma{i-1} \cup \nPi{i-1}$ and is closed under connectives, sharply bounded quantification,
  and unbounded existential (resp. universal) quantification over both the string and alphabet sorts.
\end{definition}

Thus $\nSigma{i}\subseteq \CE{i}$ and $\nPi{i}\subseteq \CA{i}$.

We denote the \emph{witnesses}, which are elements of the string sort, by Greek letters $\mu, \eta, \dots$.

\begin{definition}[Witnessing relation]
  Let $i\ge 1$.  For $A\in\CE{i}$ and a string variable $\mu$, the \emph{witnessing relation}
  $\mu\forcesS{i} A$ is defined recursively as follows:
  \begin{align*}
    \mu\forcesS{i} A
      &:= A \qquad \text{if $A\in\nSigma{i-1}\cup\nPi{i-1}$},\\
    \mu\forcesS{i} (A\land B)
      &:= \pi_0(\mu)\forcesS{i} A \land \pi_1(\mu)\forcesS{i} B,\\
    \mu\forcesS{i} (A\lor B)
      &:= (\pi_0(\mu)=0 \land \pi_1(\mu)\forcesS{i} A)
        \lor (\pi_0(\mu)=1 \land \pi_1(\mu)\forcesS{i} B),\\
    \mu\forcesS{i} \forall X\lsubseteq T\,A[X]
      &:= \forall X\lsubseteq T\,\bigl(\operatorname{nth}(X,\mu)\forcesS{i} A[X]\bigr),\\
    \mu\forcesS{i} \exists X\lsubseteq T\,A[X]
      &:= \pi_0(\mu)\lsubseteq T \land \pi_1(\mu)\forcesS{i} A[\pi_0(\mu)],\\
    \mu\forcesS{i} \exists X\,A[X]
      &:= \pi_1(\mu)\forcesS{i} A[\pi_0(\mu)],\\
    \mu\forcesS{i} \exists x\,A[x]
      &:= \pi_1(\mu)\forcesS{i} A[\operatorname{head}(\pi_0(\mu))].
  \end{align*}
  If $B\in\CA{i}$, we define \emph{co-witnessing relation} by
  \[
    \mu\forcesP{i} B := \neg(\mu\forcesS{i} \lnot B).
  \]
  We call $\mu$ appearing in $\mu\forcesS{i} A$ a \emph{witness} for $A$,
  and $\mu$ appearing in $\mu\forcesP{i} B$ a \emph{co-witness} for $B$.
\end{definition}

For a sequent $\Gamma=A_1,\ldots,A_k$ and strings $\vec\mu=\mu_1,\ldots,\mu_k$, we write
$\vec\mu\forcesS{i}\Gamma$ for $\mu_1\forcesS{i}A_1,\ldots,\mu_k\forcesS{i}A_k$, and similarly for $\forcesP{i}$.
It is easy to see that both $\mu\forcesS{i}A$ and $\mu\forcesP{i}B$ are $\nDelta{i}$-formulas over $T+\NIA{\nSigma{1}}$.

\begin{lemma}\label{lemma:casefun}
  Let $A[\vec X,\vec x]\in\CE{i}$ have all its free variables displayed.
  The function
  \[
    (\mu, \vec X, \vec x) \mapsto
      \begin{cases}
        1 & \text{if } \mu\forcesS{i}A[\vec X,\vec x],\\
        0 & \text{otherwise}
      \end{cases}
  \]
  belongs to $\Square{i}[\Alp]$ over $T+\NIA{\nSigma{1}}$.
\end{lemma}

This terminology involves a slight abuse of notation.  Given a theory $U$, a natural number $i$, and a function $F$,
we say that \emph{$F$ belongs to $\Square{i}[\Alp]$ over $U$} if
it belongs to $\Square{i}[\Alp]$ \emph{internally} to $U$ in the sense of the preceding section.
More precisely, this means that there are a $\nSigma{i-1}$-formula $\Phi$ (if $i > 1$), a code $M$ for an
oracle alphabetic Turing machine, natural numbers $n,m$, and $C,D$ such that
$\CompRes^{\Phi}_M(X^{\times n}C,X^{\times m}D,X,Y)$ is provably equivalent to the graph of $F$ over $U$.

\begin{proof}
  By the recursive definition of the witnessing relation,
  $\mu\forcesS{i}A[\vec X,\vec x]$ is obtained from
  $\nSigma{i-1}$- and $\nPi{i-1}$-formulas using Boolean connectives and sharply
  bounded quantification.  Its characteristic function is therefore computable
  with a $\ClassSigma{i-1}[\Alp]$-oracle.
\end{proof}

\begin{theorem}[Witnessing of $(\CE{i}\cup\CA{i})$-consequences]
  Let $i \ge 1$.
  Assume that $T$ is $\forall\exists!$-axiomatizable.  If a sequent
  \[
    \Sigma,\Pi
  \]
  is derivable in $T+\NIA{\nSigma{i}}$, where $\Sigma\subseteq\CE{i}$ and $\Pi\subseteq\CA{i}$,
  then there are functions $\vec f$ in $\Square{i}[\Alp]$ such that the sequent
  \[
    \vec f(\vec Z,\vec z,\vec\mu)\forcesS{i}\Sigma,
    \vec\mu\forcesP{i}\Pi
  \]
  is derivable in $T+\NIA{\nSigma{i}}$. Here $\vec Z$ and $\vec z$ list all free variables occurring in $\Sigma$ and $\Pi$.
  Moreover, if $f_s$ witnesses bounded formula $A_s\in\Sigma$, there is a polynomial bound $p_s(n)$ such that
  \[
    \size{f_s(\vec Z,\vec z,\vec\mu)} \le p_s(\size{\vec Z_{A_s}})
  \]
  is provable in $T+\NIA{\nSigma{i}}$, where $\vec Z_{A_s}$ lists all the free string variables appearing in $A_s$.
\end{theorem}
The proof is given using the standard witness argument by Buss~\cite{Buss1986}.
\begin{proof}
  By the standard free-cut elimination theorem,
  we may assume that the given derivation $D$ is anchored,
  that is, every cut formula is an instance of an axiom of either $T$ or $\NIA{\nSigma{i}}$.

  We give the representative structural and logical cases, followed by the two critical anchored-cut cases.
  The remaining cases, including the dual cases, are obtained by the same constructions.
  A complete case analysis is given in Chapter~3 of the author's master's thesis~\cite{Saito2026}.

  Hereafter, we omit displaying the free variables $\vec Z$ and $\vec z$ when they play no role.

  \paragraph{\bfseries Initial sequents and structural rules}
  For an atomic initial sequent $A,\neg A$, use empty string as the witness for $A$.
  The resulting sequent is
  \[
    \emptyseq\forcesS{i}A,\qquad \mu\forcesP{i}\neg A,
  \]
  which is derivable directly from the definitions, and polynomial bound is immediate.

  Suppose that the last inference is
  \begin{prooftree}
    \AxiomC{}
    \DeduceC{$A,A,\Gamma$}
    \UnaryInfC{$A,\Gamma$.}
  \end{prooftree}
  If $A\in\CE{i}$, the induction hypothesis for the upper sequent supplies two witnesses $g_0$ and $g_1$ for $A$.
  Define
  \[
    g:=
    \begin{cases}
      g_0,&\text{if }g_0\forcesS{i}A,\\
      g_1,&\text{otherwise}.
    \end{cases}
  \]
  By Lemma~\ref{lemma:casefun}, the decide function of $g_0\forcesS{i}A$ belongs to $\Square{i}[\Alp]$.
  Thus the function $g$ also belongs to $\Square{i}[\Alp]$ and gives a witness for $A$ in the lower sequent.
  Adding two polynomial bounds corresponding to $g_0$ and $g_1$ gives a polynomial bound for $g$.
  If $A\in\CA{i}$, simply substitute the co-witness for the lower occurrence of $A$ for both corresponding
  co-witnesses in the upper sequent.

  Weakening and exchange are immediate.

  \paragraph{\bfseries Conjunction}

  Suppose that the last inference is a conjunction:
  \begin{prooftree}
    \AxiomC{}
    \DeduceC{$\Sigma_0,\Pi_0,A$}
    \AxiomC{}
    \DeduceC{$\Sigma_1,\Pi_1,B$}
    \BinaryInfC{$\Sigma_0,\Sigma_1,\Pi_0,\Pi_1,A\land B$.}
  \end{prooftree}
  Suppose that $A\land B\in\CE{i}$.  By the induction hypothesis, there are
  $\Square{i}[\Alp]$-functions $\vec f^0,g^0$ and $\vec f^1,g^1$ such that sequents
  \[
    \begin{aligned}
      &\vec\mu_0\forcesP{i}\Pi_0,\quad
        \vec f^0(\vec\mu_0)\forcesS{i}\Sigma_0,\quad
        g^0(\vec\mu_0)\forcesS{i}A;\\
      &\vec\mu_1\forcesP{i}\Pi_1,\quad
        \vec f^1(\vec\mu_1)\forcesS{i}\Sigma_1,\quad
        g^1(\vec\mu_1)\forcesS{i}B
    \end{aligned}
  \]
  are derivable.  Define
  \[
    g(\vec\mu_0,\vec\mu_1)
      :=\operatorname{pair}\bigl(g^0(\vec\mu_0),g^1(\vec\mu_1)\bigr).
  \]
  Then we have
  \[
    \vec\mu_0\forcesP{i}\Pi_0,\quad
    \vec\mu_1\forcesP{i}\Pi_1,\quad
    \vec f^0(\vec\mu_0)\forcesS{i}\Sigma_0,\quad
    \vec f^1(\vec\mu_1)\forcesS{i}\Sigma_1,\quad
    g(\vec\mu_0,\vec\mu_1)\forcesS{i}A\land B.
  \]

  Since the size of the pair function is defined as linear in the sum of the lengths of the sizes of its components,
  the polynomial bound for $g$ can be obtained by adding the polynomial bounds for $g^0$ and $g^1$.
  
  The negative case ($A \land B \in\CA{i}$) and the dual case are analogous;
  for disjunction one attaches $0$ or $1$ to the witness of the selected disjunct.

  \paragraph{\bfseries Unbounded quantifiers}

  Suppose that the last rule introduces a string existential:
  \begin{prooftree}
    \AxiomC{}
    \DeduceC{$\Sigma,\Pi,A[U]$}
    \UnaryInfC{$\Sigma,\Pi,\exists X\,A[X]$.}
  \end{prooftree}
  Then we have $A[X] \in \CE{i}$. If the induction hypothesis gives a witness $g(\vec\mu)$ for $A[U]$, then
  $h(\vec{\mu}) := \operatorname{pair}(U,g(\vec\mu))$ witnesses $\exists X\,A[X]$, if the principal formula is unbounded.
  If bounded, say, $\exists X\,A[X] = \exists X\le T\,A'[X]$,
  then witness is given by
  \[
     h(\vec{\mu}) := \begin{cases}
       \operatorname{pair}(U, g(\vec\mu)) & \text{if $U \le T$},\\
       \emptyseq & \text{otherwise.}
     \end{cases}.
  \]
  The polynomial bounds in both cases are given by $O(\size{T} + p_g(\size{T} + \size{\vec{Z}_{\exists X\,A[X]}}))$,
  where $p_g$ is the polynomial bound for the witness $g$.

  For an alphabet existential introduced with a term $t$, we use $\operatorname{pair}(\brak{t},g(\vec\mu))$ instead.
  
  Universal quantifiers are handled similarly, using the appropriate projections of the witness for the upper sequent.

  \paragraph{\bfseries Sharply bounded quantification}

  Suppose that the last rule introduces a sharply bounded universal quantifier:
  \begin{prooftree}
    \AxiomC{}
    \DeduceC{$\Sigma,\Pi,Y\not\lsubseteq T,A[Y]$}
    \UnaryInfC{$\Sigma,\Pi,\forall X\lsubseteq T\,A[X]$.}
  \end{prooftree}
  If $A \in \CE{i}$, then the induction hypothesis gives functions $\vec{f}, g$ such that the sequent
  \[
    \vec{f}(\vec{\mu}, Y)\forcesS{i}\Sigma,\quad
    \vec{\mu}\forcesP{i}\Pi,\quad
    Y \lsubseteq T \to g(\vec{\mu}, Y)\forcesS{i}A[Y]
  \]
  is derivable.
  Form a witness $h(\vec\mu)$ by coding the sequence $\langle g(\vec{\mu}, I) \rangle_{I \lsubseteq T}$
  and using it as the witness for $\forall X\lsubseteq T\,A[X]$.
  Define $f'_s(\vec{\mu})$ as $f_s(\vec{\mu}, I)$ for the first $I \lsubseteq T$ such that $\lnot g(\vec{\mu}, I) \forcesS{i} A[I]$,
  otherwise use empty strings. Then sequent
  \[
    \vec{f}'(\vec{\mu})\forcesS{i}\Sigma,\quad
    \vec{\mu}\forcesP{i}\Pi,\quad
    h(\vec{\mu})\forcesS{i}\forall X\lsubseteq T\,A[X]
  \]
  is derivable.
  The polynomial bound for $h(\vec{\mu})$ is obtained by $O({1\over 2}\size{T} (\size{T} + 1) p_{g}(\size{\vec{Z}_{\forall X\lsubseteq T\,A[X]}}))$,
  where $p_{g}$ is the polynomial bound for the witness $g$.
  The sharply bounded existential and the negative sharply bounded cases are analogous.

  Universal axioms of $T$ and the universal basic axioms for strings are
  handled as initial sequents, since they require no witness.  It remains to
  consider cuts anchored on the existential axioms of $T$ and on instances of
  induction.
  
  \paragraph{\bfseries Case when the last rule is a cut anchored on an axiom of $T$}
  
  Assume that $D$ is of the form
  \begin{prooftree}
    \AxiomC{}
    \DeduceC{$\Sigma, \Pi, \forall\vec{y}\, \lnot A[\vec{t}, \vec{y}]$}
    \AxiomC{}
    \RightLabel{\scriptsize $\mathrm{ax}$}
    \UnaryInfC{$\exists\vec{y}\, A[\vec{t}, \vec{y}]$}
    \RightLabel{\scriptsize cut}
    \BinaryInfC{$\Sigma,\Pi$}
  \end{prooftree}
  where $A$ is quantifier-free and $T$ proves that $\vec y$ is unique.
  Also notice that $\forall\vec y\, \lnot A(\vec t, \vec y) \in \CA{i}$.
  By the induction hypothesis, there exist functions $\vec f$ and $g$ in $\Square{i}[\Alp]$ such that
  \[
    \vec\mu\forcesP{i}\Pi,\quad
    \vec f(\vec Z, \vec z, \vec\mu,\sigma)\forcesS{i}\Sigma,\quad
    \sigma\forcesP{i}\forall\vec y\,\lnot A[\vec t,\vec y]
  \]
  is derivable.
  Since $T$ proves the uniqueness of $\vec y$,
  we have a unit-time $\Alp$-alphabetic TM $M$ that computes the unique $\vec y$ such that $A[\vec t,\vec y]$ holds,
  given $\vec{z}$ as input.
  Thus we have $M(\vec{z}) \forcesS{i} \exists\vec{y}\, A[\vec{t},\vec{y}]$,
  and the following sequent is derivable by a cut.
  \[
    \vec\mu\forcesP{i}\Pi,\quad
    \vec g(\vec Z,\vec z, \vec\mu) \forcesS{i}\Sigma.
  \]
  Here $\vec g(\vec Z,\vec z, \vec\mu) := \vec f(\vec Z,\vec z, \vec\mu,M(\vec z))$.
  The polynomial bound for $\vec g$ is immediate.

  \paragraph{\bfseries Case when the last rule is a cut anchored on an axiom of $\NIA{\nSigma{i}}$}

  Assume that $D$ is of the form
  \begin{prooftree}
    \AxiomC{}
    \RightLabel{$D'$}
    \DeduceC{$\Sigma,\Pi,A[\emptyseq] \land \forall Y\forall y\,(A[Y] \to A[Y y]) \land \lnot A[T]$}
    \AxiomC{}
    \RightLabel{\scriptsize $\mathrm{ax}$}
    \UnaryInfC{$\Phi$}
    \RightLabel{\scriptsize cut}
    \BinaryInfC{$\Sigma,\Pi$}
  \end{prooftree}
  Here $A \in \nSigma{i}$ and $T$ is a string term.

  By applying $\forall$, $\lor$, and $\land$ inversions to $D'$,
  we have anchored derivations of sequents
  \begin{gather*}
    \Pi,\Sigma,A[\emptyseq];\\
    \Pi,\Sigma,\neg A[Y],A[Yy];\\
    \Pi,\Sigma,\neg A[T]
  \end{gather*}
  where $Y$ and $y$ are fresh, with smaller heights than $D$.
  The induction hypothesis therefore
  supplies $\Square{i}[\Alp]$-functions for which the following sequents
  are derivable:
  
  \begin{align*}
    \Iota :=&\quad \vec\mu\forcesP{i}\Pi,\quad
      \vec f^0(\vec\mu)\forcesS{i}\Sigma,\quad
      g^0(\vec\mu)\forcesS{i}A[\emptyseq]\\
    \Lambda :=&\quad \vec\mu\forcesP{i}\Pi,\quad
      \vec f^s(Y,y,\vec\mu,\nu)\forcesS{i}\Sigma,\quad
      \nu\forcesP{i}\neg A[Y],\quad
      g^s(Y,y,\vec\mu,\nu)\forcesS{i}A[Yy]\\
    \Gamma :=&\quad \vec\mu\forcesP{i}\Pi,\quad
      \vec f^t(\vec\mu,\nu)\forcesS{i}\Sigma,\quad
      \nu\forcesP{i}\neg A[T]
  \end{align*}

  Define $g$, $\vec f$, and $\vec h$ as follows:
  \begin{align*}
    g(\vec\mu,Y)
      &:=
      \begin{cases}
        g^0(\vec\mu),
          &\text{if }Y=\emptyseq,\\
        g^s\bigl(Y',y,\vec\mu,g(\vec\mu,Y')\bigr),
          &\text{if }Y=Y'y,
      \end{cases}\\
    f_j(\vec\mu,Y)
      &:=
      \begin{cases}
        f^0_j(\vec\mu),
          &\text{if }Y=\emptyseq,\\
        f^s_j\bigl(Y',y,\vec\mu,g(\vec\mu,Y')\bigr),
          &\text{if }Y=Y'y
            \text{ and }g(\vec\mu,Y')\forcesS{i}A[Y'],\\
        f_j(\vec\mu,Y'),
          &\text{if }Y=Y'y
            \text{ and }\neg\bigl(g(\vec\mu,Y')\forcesS{i}A[Y']\bigr),
      \end{cases}\\
    h_j(\vec\mu)
      &:=
      \begin{cases}
        f^t_j\bigl(\vec\mu,g(\vec\mu,T)\bigr),
          &\text{if }g(\vec\mu,T)\forcesS{i}A[T],\\
        f_j(\vec\mu,T),
          &\text{if }\neg\bigl(g(\vec\mu,T)\forcesS{i}A[T]\bigr).
      \end{cases}
  \end{align*}
  Since $A\in\nSigma{i}$, the induction hypothesis gives polynomial bounds $p_0,p_s$ such that,
  \[
    \size{g^0(\vec\mu)}\le p_0(\size{\vec Z_{A[\emptyseq]}}),\qquad
    \size{g^s(Y,y,\vec\mu,\nu)}\le p_s(\size{\vec Z_{A[\emptyseq]}}+\size{Y}).
  \]
  Since the polynomial bound of $g^s$ does not depend on $\nu$, The function $g$ can be defined by the bounded recursion.
  The polynomial bound for $g$ can be taken as
  \[
    p_0(\size{\vec Z_{A[\emptyseq]}}) + p_s(\size{\vec Z_{A[\emptyseq]}} + \size{Y}).
  \]
  The polynomial bounds for the functions $\vec f$ and $\vec h$ can be obtained immediately.

  We first show, by induction on notation, that the sequent
  \[
    \Delta :=
      \vec\mu\forcesP{i}\Pi,\quad
      \bigvee \vec f(\vec\mu,Y)\forcesS{i}\Sigma
      \lor g(\vec\mu,Y)\forcesS{i}A[Y]
  \]
  is derivable.
  The case $Y = \emptyseq$ is precisely $\Iota$.
  Suppose that $Y = Y'y$.
  The two cases in the definition of $\vec f$ give the following two sequents:
  \begin{align*}
    \Delta_0 := {}&
      \vec\mu\forcesP{i}\Pi,\quad
      g(\vec\mu,Y')\forcesP{i}\neg A[Y'],\\
    &\Bigl[
      \bigvee\vec f(\vec\mu,Y')\forcesS{i}\Sigma
      \lor g(\vec\mu,Y')\forcesS{i}A[Y']
      \Bigr]\\
    {}\to{}&
      \Bigl[
      \bigvee\vec f^s(Y',y,\vec\mu,g(\vec\mu,Y'))\forcesS{i}\Sigma
      \lor g^s(Y',y,\vec\mu,g(\vec\mu,Y'))\forcesS{i}A[Y'y]
      \Bigr];\\
    \Delta_1 := {}&
      \vec\mu\forcesP{i}\Pi,\quad
      g(\vec\mu,Y')\forcesS{i}A[Y'],\\
    &\Bigl[
      \bigvee\vec f(\vec\mu,Y')\forcesS{i}\Sigma
      \lor g(\vec\mu,Y')\forcesS{i}A[Y']
      \Bigr]\\
    {}\to{}&
      \Bigl[
      \bigvee\vec f(\vec\mu,Y')\forcesS{i}\Sigma
      \lor g^s(Y',y,\vec\mu,g(\vec\mu,Y'))\forcesS{i}A[Y'y]
      \Bigr].
  \end{align*}
  The sequent $\Delta_0$ follows from $\Lambda$ by substituting
  $g(\vec\mu,Y')$ for $\nu$, and $\Delta_1$ is a tautology. This proves $\Delta$.
  Moreover, by the definition of $\vec h$, we have
  \[
    \begin{aligned}
      \Xi_0 :=&\quad
        \vec\mu\forcesP{i}\Pi,\quad
        g(\vec\mu,T)\forcesP{i}\neg A[T],\quad
        \vec f^t(\vec\mu,g(\vec\mu,T))\forcesS{i}\Sigma,\\
      \Xi_1 :=&\quad
        \vec\mu\forcesP{i}\Pi,\quad
        g(\vec\mu,T)\forcesS{i}A[T],\quad
        \vec f(\vec\mu,T)\forcesS{i}\Sigma.
    \end{aligned}
  \]
  The sequent $\Xi_0$ follows from $\Gamma$ by substituting
  $g(\vec\mu,T)$ for $\nu$, and $\Xi_1$ follows from $\Delta[Y := T]$.
  A final case distinction therefore yields
  \[
    \vec\mu\forcesP{i}\Pi,\qquad
    \vec h(\vec\mu)\forcesS{i}\Sigma
  \]
  as required.
\end{proof}

We assume that $i\ge 1$.

\begin{corollary}[Witnessing for $T+\NIA{\nSigma{i}}$]
  Let $A[X,Y]$ be a $\nSigma{i-1}$-formula.  If
  \[
    T+\NIA{\nSigma{i}} \vdash \forall X\,\exists Y\,A[X,Y],
  \]
  then there is a function $f\in\Square{i}[\Alp]$ such that
  \[
    T+\NIA{\nSigma{i}} \vdash \forall X\,A[X,f(X)].
  \]
\end{corollary}

\begin{corollary}
  Assume that the theory $T$ is $\forall\exists!$-axiomatizable and $\forall\exists$-complete for $\mathbb{K}$.
  The provably total $\nSigma{i}$-functions in $T+\NIA{\nSigma{i}}$ are exactly the functions in $\Square{i}[\mathbb{K}]$.
\end{corollary}

Specifically,

\begin{corollary}
  \leavevmode
  \begin{enumerate}
    \item The provably total $\nSigma{i}$-functions in $\BOOL+\NIA{\nSigma{i}}$ are exactly the functions in $\Square{i}$.
    \item The provably total $\nSigma{i}$-functions in $\RCOF+\NIA{\nSigma{i}}$ are exactly the functions in $\Square{i}[\Real]$.
  \end{enumerate}
\end{corollary}

\section*{Acknowledgements}
I would like to thank my supervisor Professor Keita Yokoyama and Professor Takashi Yamazaki,
for their guidance and support throughout this research.
I also thank Professor Satoru Kuroda for taking the time to discuss my research and for suggesting relevant topics. 
The author used AI to review and improve its English grammar and writing.

\bibliographystyle{plainnat}
\bibliography{references}

\end{document}